\documentclass[10pt]{article}
\usepackage{latexsym,amssymb}

\def\N{\mathbb{N}}
\def\Q{\mathbb{Q}}
\def\Z{\mathbb{Z}}
\def\R{\mathbb{R}}
\def\T{\mathbb{T}}
\def\D{\mathbb{D}}
\def\P{\mathbb{P}}
\def\fl{\mathrm{flim}}
\def\hl{\mathrm{hlim}}
\def\proof{\par\medskip\noindent{\em Proof. }}
\def\eproof{\hfill{$\Box$}\bigskip}

\newtheorem{thm}{Theorem}[section]
\newtheorem{pro}[thm]{Proposition}
\newtheorem{lem}[thm]{Lemma}
\newtheorem{ddef}[thm]{Definition}

\begin{document}

\begin{center}
{\bf\large Real numbers as infinite decimals and irrationality of $\sqrt{2}$}

\bigskip\
{\large Martin Klazar}\\

\medskip
{\em
Department of Applied Mathematics (KAM)\\
and\\ 
Institute for Theoretical Computer Science (ITI)\\
Charles University, Faculty of Mathematics and Physics\\
Malostransk\'e n\'am\v est\'\i\ 25, 118 00 Praha\\
Czech Republic\\
{\tt klazar at kam.mff.cuni.cz}
}
\end{center}

\bigskip
\begin{abstract}
In order to prove irrationality of $\sqrt{2}$ by using only decimal expansions (and not fractions), we 
develop in detail a model of real numbers based on infinite decimals and arithmetic operations with them. 
\end{abstract}

\section{Introduction} 

{\bf Irrationality of $\sqrt{2}$ via decimals?}
Many proofs of irrationality of the number $\sqrt{2}$ were published and collected, e.g.,
Beigel \cite{beig}, Bogomolny \cite{zwebu} (at least 19 on-line proofs), Flannery \cite{flan}, 
Gardner \cite{gard} (lists 18 references), 
Harris \cite{harr} (13 proofs), Miller and Montague \cite{mill_mont}, Myerson \cite{myer}, Subbarao \cite{subb}, and
Waterhouse \cite{wate}, and many textbooks on algebra, mathematical analysis and number theory 
present such proof: Allouche and Shallit \cite[Theorem 2.2.1]{allo_shal}, Apostol \cite[Theorem 1.10]{apos}, 
Hardy and Wright \cite[Theorem 43]{hard_wrig}, Jarn\'\i k \cite[str. 16]{jarn}, Pugh \cite[p. 11]{pugh}, Tao \cite[Proposition 4.4.4]{tao_ma}, 
Zorich \cite[p. 51]{zori}, to name a few. These proofs start invariably from the hypothetical 
fractional representation $\sqrt{2}=\frac{a}{b}$ with $a,b\in\Z$ and derive a contradiction. But could one prove that 
$$
\sqrt{2}=1.414213562373095048801688724209698\dots
$$ 
is irrational in a purely numerical way by showing that this decimal expansion is not ultimately periodic, without 
invoking the fractional model of rational numbers? Such proofs are hard to come 
by in the literature. We know of only two imperfect examples and reproduce one of them below. (The other is proof no. 
12 in \cite{zwebu}, due to A. Cooper, which is less convincing as it fails to deal, without help of fractions, with the 
case 2 below.)

\medskip\noindent
{\bf A proof of irrationality of $\sqrt{2}$ by decimals.}
Lindstrom \cite{lind} found the following proof of irrationality of $\sqrt{2}$ using 
decimal expansions. It brings to contradiction the assumption that $\sqrt{2}$ 
is a terminating or repeating decimal. We abridge it somewhat but preserve its style and notation.
\begin{quote}
Case 1. If $\sqrt{2}$ is a terminating decimal, say $\sqrt{2}=1.abc$ with decimal digits
$a,b,c$ and $c$ is nonzero, multiplying by $1000$ and squaring, we get ($1abc$ is in the 
decimal notation)
$$
2000000=(1000\sqrt{2})^2=(1abc)^2.
$$
This cannot hold as the units digit of $2000000$ is zero but the units digit of the right-hand 
side is nonzero, the units digit of $c^2$. Other terminating decimals are treated 
in a similar way.

\noindent
Case 2. If $\sqrt{2}$ is a repeating decimal, say 
$\sqrt{2}=1.\overline{abc}=1.abcabc\dots$ with some decimal digits $a,b,c$, then  ($abc$ is in the decimal notation)
$\sqrt{2}=1+abc/10^3+abc/10^6+\dots\;$. Summing the geometric series, rearranging and 
squaring, we get
$$
2\cdot999^2=(\sqrt{2}(10^3-1))^2=((10^3-1)+abc)^2.
$$
This cannot hold either as the units digit on the left-hand side after multiplying is $2$ but the units digit on
the right-hand side after squaring is that of a square, one of $0,1,4,5,6,9$. Other repeating decimals 
are treated in a similar way.

\noindent
Hence $\sqrt{2}$ is irrational.
\end{quote}
{\bf Is anything wrong with the proof?} It is not and it is. In the liberal view that allows mixing infinite decimals 
with fractions, the above proof is fine. Indeed, to be fair, it appears that its author did not set as his goal to 
use exclusively decimal expansions. But in the more restricted view when one accepts only arguments using decimal 
expansions, the above proof suffers from two problems. First, in summing 
$1+abc/10^3+abc/10^6+\dots$ it uses fractions. Second, and this is a more serious problem, it uses 
the familiar but in our context problematic simplification $(\sqrt{2})^2=2$. When real numbers are infinite decimals, 
the correct simplification is
$(\sqrt{2})^2=\{1.999\dots,2\}$, i.e., one has to consider two possibilities, $(\sqrt{2})^2=1.999\dots$ and 
$(\sqrt{2})^2=2$. Indeed, $(1.4142135\dots)^2$, calculated via truncations as $1.4^2$, $1.41^2$, $1.414^2$ and so 
on, approaches $1.999\dots$ and not $2$. The above proof does not consider the possibility 
$(\sqrt{2})^2=1.999\dots$ at all and therefore is incomplete. 

But wait, may object the attentive reader, isn't it so that we identify $1.999\dots$ with $2$ and thus it suffices to 
work just with the representative $2$? Therefore, isn't the argument in case 2 looking only at 
$(\sqrt{2})^2=2$ sufficient after all? It is and it is not. It is sufficient because when 
the model of $\R$ based on decimal expansions is worked out in detail (in Section 2), it turns out that it 
really suffices to look, in a sense, only at the possibility $(\sqrt{2})^2=2$ (see the proofs of 
Propositions~\ref{odm2vr} and \ref{odm2vd}). 
At the same time it is not sufficient because if such model is not specified or at least referred to, which, alas, is 
the rule in discussions on this topic, all arguments about what 
can or cannot be justified by identifications like $1.999\dots=2$ are necessarily scholastic and lack substance 
(cf. various ``proofs'' that $0.999\dots=1$ at \cite{wiki09je1}).

\medskip\noindent
{\bf What is in this text.} Our original motivation was to understand if the proof 
of Lindstrom qualifies as a ``numerical'' proof of irrationality of $\sqrt{2}$ using only decimal expansions
and if it does not, how would such a proof look like. As to Lindstrom's proof, we summarized our opinion in the two preceding paragraphs. In Section 3 we present our proofs of 
irrationality of $\sqrt{2}$ using only decimal expansions. They require a model of $\R$ based on decimal expansions, 
called {\em a decimal model of $\R$} for short, which we develop in detail in Section 2.

Such model is of course not a complete novum because besides the two well known 1872 models of $\R$, the 
model of M\'eray, Cantor and Heine using fundamental sequences of rationals and 
Dedekind's model based on cuts on rational numbers, there is, allegedly, a third model 
based on decimal expansions, developed by Weierstrass and Stolz. Unfortunately, unlike for the first two classical 
models, we could not find in the literature any really detailed and satisfactory presentation of the Weierstrass--Stolz model or its modern version, only outlines and sketches; we give references and more comments at the end of Section 2 where
we present our stab at this model.

\medskip\noindent
{\bf Troubles with formal limits.} So how to add and multiply infinite decimals? Fowler \cite{fowl} challenges the 
reader to calculate  the product
$$
1.222222\dots\times0.818181\dots=\;?
$$
A natural idea is to take finite truncations of the factors, form the 
sequence of their products,
$$
1.2\times0.8=0.96,\ 1.22\times0.81=0.9882,\ 1.222\times0.818=0.999596,\ \dots,
$$
and take its formal limit. Each digit after the decimal point stabilizes 
eventually at $9$ and each digit before it is always $0$. Hence, indisputably, 
$$
1.222222\dots\times0.818181\dots=0.999999\dots\;.
$$
In this way, using formal convergence of digits, one can define multiplication and addition of infinite decimals.  

This often suggested approach brings certain troubles. It takes an effort to show that 
subtraction is well defined because the sequence of differences of the truncations is in general not eventually monotone. 
But the main trouble is that formal limits do not commute with arithmetic operations and consequently these operations
on infinite decimals lose their convenient properties.  The following examples show that addition ceases to be 
associative and that the distributive law fails:
$$
(-0.999\dots+1)+0.999\dots=0.999\dots\ \mbox{ but }\ -0.999\dots+(1+0.999\dots)=1
$$ 
and 
$$
(10-1)\cdot0.111\dots=0.999\dots\ \mbox{ but }\ 10\cdot0.111\dots-1\cdot0.111\dots=1.
$$
Certainly, when identifications like $0.999\dots=1$ are applied these irregularities go away (as one expects they must) 
but again it takes certain effort to prove it rigorously (without handwaving). 

We compare commutativity and associativity of arithmetic operations on infinite decimals. Since
these operations on the finite truncations are clearly commutative, applying formal limits we get that 
addition and multiplication of infinite decimals is commutative. At first one might think that, similarly, so it
works for associativity (cf. Gowers \cite{gowe_dec}) but, with the previous example, we know very
well that this argument is erroneous. Addition of finite decimal truncations is of course associative but 
formal limits do not transfer this to infinite decimals. The point is that (unlike 
commutativity) associativity of an operation with infinite decimals amounts to the exchangeability 
of two formal limits, which is a nontrivial
result to be proven, and  consequently is not at all automatic from the associativity for finite truncations. For multiplication this exchangeability occurs and multiplication of infinite decimals is 
associative. For addition it fails and addition is not associative; the same happens for the distributive law.
(The reason is that for multiplication one has always monotonicity of partial products of truncations, which is in 
general not the case for subtraction.)

Let us also illustrate the difference between usual metric limit and formal limit of decimal expansions by a 
curious paragraph from Courant and Robbins \cite[p. 293]{cour_robb}. After explaining the standard notion of the limit 
$a$ of a real sequence $a_1,a_2,\dots,$ they write
\begin{quote}
If the members of the sequence $a_1,a_2,a_3,\cdots$ are expressed as infinite decimals, then the statement 
$\lim a_n=a$ simply means that for any positive integer $m$ the first $m$ digits of $a_n$ coincide with the first 
$m$ digits of the infinite decimal expansion of the fixed number $a$, provided that $n$ is chosen sufficiently 
large, say greater than or equal to some value $N$ (depending on $m$). This merely corresponds to choices of $\epsilon$ 
in the form $10^{-m}$.
\end{quote}
So, since the first $m$ digits of $a_n$ are for large $n$ equal to the first $m$ digits of {\em the} decimal expansion 
of $a$, should one conclude that each digit of $a_n$ eventually stabilizes and the decimal expansions of 
$a_n$ formally converge? It is not exactly what the paragraph says (though it may seem to follow immediately from it)
and it is not actually true, as shows the sequence $a_1=1.1, a_2=0.9, a_3=1.01, a_4=0.99, a_5=1.001, a_6=0.999, \dots$ 
that metrically converges to $a=1$ but has no formal limit (replacements of $0.999$ by $0.998999\dots$ or $1.001$ by 
$1.000999\dots$ etc. change nothing on this). The correct reading is that the first $m$ digits of $a_n$ are for large 
$n$ equal to the first $m$ digits of {\em a} decimal expansion of $a$. In our example, one needs to switch constantly 
between the two expansions $0.999\dots$ and $1.000\dots$ of $a=1$, the former being used for $a_{2n}$ and the latter for $a_{2n-1}$. Is this what the authors meant? (No further explanation or example are given.)

\medskip\noindent
{\bf Our approach.} We had tried formal limits as well but then abandoned them in favor of an alternative
approach. Instead of developing the cumbersome arithmetics of infinite decimals and fixing it at the end by
identifications like $0.999\dots=1$, we work with them from the very beginning and define arithmetic
operations with infinite decimals by means of a limit weaker than the formal limit as inherently multivalued
(i.e., sometimes bivalued) operations. Our answer to Fowler's question is that
$$
1.222222\dots\times0.818181\dots=\{0.999999\dots,1\}
$$
and our definition of multiplication (and addition) is such that both possibilities for the result are on completely 
equal footing. Similarly we define, for example, that 
$$
\{0.1999\dots,0.2\}+\{-0.5,-0.4999\dots\}=\{-0.3,-0.2999\dots\}.
$$
This is a less ambitious approach than formal limits as we do not care which of the two
possibilities of the result is the ``correct'' one. Its advantage is that all 
required properties of arithmetic operations with real numbers can be established in a straightforward and  
natural way and thus, finally, a sound numerical proof of the irrationality of $\sqrt{2}$ can be formulated. 
However, at the end of Section 2 we return to formal limits and prove that they do provide well defined arithmetic operations with infinite decimals.

\section{A decimal model of $\R$} 

We start from the ordered ring of integers $\Z=(\Z,+,\cdot,<)$, which we assume to be given, and build from it
the complete ordered field of real numbers $\R=(\R,+,\cdot,<)$ and, on the way, the field of rational numbers 
$\Q=(\Q,+,\cdot)$. By a ring we always mean a commutative ring with $1$. We use notation $\N=\{1,2,\dots\}$ 
(natural numbers) and $\Z=\{\dots,-2,-1,0,1,2,\dots\}$ (the integers).

\medskip\noindent
{\bf Decimals and real numbers.} A {\em signed decimal} $d$ is a pair of the sign $+$ or $-$ and an infinite string of decadic {\em digits} which are indexed by the integers $k,k-1,k-2,\dots$ and the first of which is nonzero:
$$
d=\pm a_ka_{k-1}a_{k-2}\dots
$$ 
where $k\in\Z$, $a_i\in\{0,1,2,3,4,5,6,7,8,9\}\subset\Z$ and $a_k\ne0$. The {\em zero decimal} $0$ has no sign, $k=0$, and 
has only zero digits:
$$
0=a_0a_{-1}a_{-2}\dots=000\dots\;.
$$
The set of all {\em decimals}, signed and zero, is denoted $\D$. Decimals with the $+$ sign are {\em positive} and those 
with the $-$ sign {\em negative}. For $d\in\D$ and $i\in\Z$ we write $d_i$ for the 
$i$th digit $a_i$ of $d$. If $a_i$ does not exist ($i>k$ for signed decimals or $i>0$ for the zero decimal) we set $d_i=0$.
The equality of two decimals $d$ and $e$ means that they have equal signs and $d_i=e_i$ for every $i\in\Z$, or both are 
zero. Note that we regard digits in decimals as integers (not just labels) and add, multiply and 
compare them as elements of $\Z$. The decimal differing from a decimal $d$ only in sign is denoted $-d$ 
(for the zero decimal we set $-0=0$). When 
writing down decimals we use usual conventions concerning the decimal point and omitting the $+$ sign and 
the trailing zeros. 

We define a linear order $(\D,<)$ on decimals. It is the lexicographic ordering, adjusted for signs. More precisely,
for $d,e\in\D$, $d\ne e$, we set $d<e$ if and only if $d$ and $e$ have the same $+$ (respectively $-$) sign and there
is an $m\in\Z$ such that $d_i=e_i$ for $i>m$ but $d_m<e_m$ (respectively $d_m>e_m$), or $d$ has $-$ sign and $e=0$ or 
$e$ has $+$ sign, or $d=0$ and $e$ has $+$ sign.

This linear ordering is not dense, for some pairs of decimals $d<e$ there is no decimal $f$ with $d<f<e$. 
It happens exactly when $d$ and $e$ have the same $+$ (respectively $-$) sign and there is an $m\in\Z$
such that $d_i=e_i$ for $i>m$, $e_m-d_m=1$ (respectively $e_m-d_m=-1$), and for $i<m$ one has $d_i=9$ and $e_i=0$ 
(respectively $d_i=0$ and $e_i=9$). We call such pairs $d<e$ {\em jumps}. Jumps are clearly disjoint. Examples of 
jumps are
$$
0.999\dots<1\ \mbox{ or }\ -17.341<-17.340999\dots\;.
$$
We write $\T_9$ 
for the decimals ending with infinitely many $9$'s and $\T_0$ for the decimals ending with infinitely many $0$'s. 
We call the latter {\em terminating decimals}. Note that $0\in\T_0$. Every jump consists of one element from 
$\T_0$ and one element from $\T_9$ and every element of $\T_0\cup\T_9\backslash\{0\}$ appears in exactly one jump.

Let $\sim$ be the equivalence relation on $\D$ identifying elements in jumps, 
i.e., $d\sim e$ iff $d=e$ or if $d<e$ or $e<d$ is a jump. We define the set of {\em real numbers} as the set 
of equivalence classes. 

\begin{ddef}
The set of real numbers $\R$ is
$$
\R=\D/\!\sim\;=\{\{0\}\}\cup\{\{d\}\;|\;d\in\D\backslash(\T_0\cup\T_9)\}\cup\{\{c,d\}\;|\;c<d\mbox{ is a jump}\}.
$$
\end{ddef}

\noindent 
For $X\subset\D$ we let $[X]$ denote the set of elements of $\R$ intersecting $X$; for $d\in\D$ we write $[d]$ 
(and not $[\{d\}]$) for the equivalence class of $d$. 
For $\alpha,\beta\in\R$ we set $\alpha<\beta$ iff $d<e$ for some decimals $d\in\alpha$ and $e\in\beta$. It is clear 
that $(\R,<)$ is a dense linear ordering.

We shall prove the following two theorems.

\begin{thm}\label{oR}
The structure $\R=(\R,+,\cdot,<)$, with the operations $+$ and $\cdot$ on the real numbers still to 
be defined, is a complete ordered field. 
\end{thm}
We let $\P$ denote the set of ultimately periodic decimals: $d\in\P$ iff there are $m\in\Z$ and $p\in\N$ such 
that $d_i=d_{i-p}$ for every $i<m$. We have $\T_0,\T_9\subset\P$.
\begin{thm}\label{oQ}
The prime field of the field $\R=(\R,+,\cdot)$, isomorphic to $\Q$, is formed exactly by the ultimately periodic 
real numbers $[\P]$.
\end{thm}

Before we introduce arithmetics of real numbers, we dispose with the completeness in
Theorem~\ref{oR} and show that in $(\R,<)$ every nonempty set bounded from above has supremum, i.e.,
the least upper bound. It is clear from how $(\R,<)$ arises from $(\D,<)$ that it suffices to prove this for  
decimals.

\begin{pro}\label{tlup}
In $(\D,<)$ every nonempty set bounded from above has the least upper bound.
\end{pro}
\proof
Let $D\subset\D$, $D\ne\emptyset$, be bounded from above. We may assume that all elements in $D$ are positive; other 
cases easily reduce to this. Since $D$ has an upper bound, there is a $k\in\Z$ such that $d_i=0$ for 
every $d\in D$ and $i>k$ but $d_k>0$ for some $d\in D$. We set $c_k=\max d_k$, taken over all $d\in D$. Hence $c_k>0$. 
We set $c_{k-1}=\max d_{k-1}$, taken over all $d\in D$ with
$d_k=c_k$. We set $c_{k-2}=\max d_{k-2}$, taken over all $d\in D$ with $d_k=c_k,d_{k-1}=c_{k-1}$. Continuing this way,
we define a decimal $c\in\D$, with $+$ sign. It is immediate that $c$ is the least upper bound of $D$. 
\eproof

\noindent
Similarly for infima, i.e., largest lower bounds.

\medskip\noindent
{\bf Arithmetic operations with real numbers.} We get them as extensions of the unproblematic arithmetic operations 
with terminating decimals, which we begin with. Let $\Z_0$ be the set of all decimals $d$ satisfying $d_i=0$ 
for every $i<0$. The mapping (decadic notation) $0\mapsto 0$ and
$$
\pm a_ka_{k-1}\dots a_0000\dots\mapsto \pm(a_k10^k+a_{k-1}10^{k-1}+\cdots+a_010^0)
$$
is a 1-1 correspondence between $\Z_0$ and $\Z$. We get an ordered ring $(\Z_0,+,\cdot,<)$, isomorphic to 
the ring of integers. For $k\in\Z$ we let $10^k$ denote the positive terminating decimal with 
$$
\mbox{$(10^k)_k=1$ and $(10^k)_i=0$ for $i\ne k$.}
$$ 
For $d\in\D$ and $k\in\Z$ we denote by $10^kd$ the decimal obtained from $d$ by keeping its sign and shifting its digits 
by $k$ places: 
$$
(10^kd)_i=d_{i-k}.
$$
Note that $10^k(10^{-k}d)=d$ for every $k\in\Z$ and $d\in\D$ and that the elements $d$ and $10^kd$ of $\Z_0$ are 
mapped by decadic notation to the elements $a$ and $10^ka$ of $\Z$.
For every $d\in\T_0$ there is an $N\in\N$ such that $10^nd\in\Z_0$ for every $n\ge N$. To add and multiply two terminating decimals $c$ and $d$, 
we take some $k\in\Z$ such that $10^kc$ and $10^kd$ lie in $\Z_0$ and define 
$$
c+d=10^{-k}(10^kc+10^kd)\ \mbox{ and }\ c\cdot d=10^{-2k}(10^kc\cdot10^kd) 
$$
where the operations on the right sides are in $\Z_0$. It follows that the results do not depend on $k$. Also, 
for $c,d\in\T_0$ we have $c<d$ iff the terminating decimal $d-c$ has $+$ sign (i.e., the order on $\Z_0$ transferred from
$\Z$ coincides with the already defined lexicographic ordering). 

\begin{pro}\label{aotd}
$(\T_0,+,\cdot,<)$ is an ordered ring.
\end{pro}
\proof
The elements $0$ and $1$ ($=1.000\dots$) are neutral to addition and multiplication and any $d\in\T_0$ has
additive inverse $-d$. The required properties of $+,\cdot$ and $<$ (commutativity, associativity, 
distributivity and monotonicity) follow from the fact that they hold in $\Z_0$. For example, if 
$c(d+e)\ne cd+ce$ for some $c,d,e\in\T_0$, we take a $k\in\N$ so that all $10^kc,10^kd$ and $10^ke$ lie in $\Z_0$ and
obtain a refutation of the distributive law in $\Z_0$:  
$10^{2k}(c(d+e))=(10^kc)((10^kd)+(10^ke))$ differs from $10^{2k}(cd+ce)=(10^kc)(10^kd)+(10^kc)(10^ke)$.
Similarly for other properties.
\eproof

\noindent
Two defects of this ring, non-completeness and the lack of division, will be fixed by extension of the operations
to $\R$.

For $n\in\N$ and $d\in\D$, the $n$-{\em truncation} $d|n\in\T_0$ has the same sign as $d$ and digits 
$$
\mbox{$(d|n)_i=d_i$ for $i\ge-n$ and $(d|n)_i=0$ for $i<-n$}.
$$ 
If $(d|n)_i=0$ for every $i\in\Z$, we omit the sign and set $d|n$ to be the zero decimal. Any truncation of the zero 
decimal is the zero decimal. For a sequence of terminating decimals $c(n)$, $n=1,2,\dots,$ we write 
$$
c(n)\to0 
$$
if for every $k\in\N$ there is an $N\in\N$ such that $c(n)|k=0$ whenever $n>N$. It is clear that $c(n)\to0$ and 
$d(n)\to0$ implies $c(n)+d(n)\to0$. Similarly, $c(n)-c'(n)\to0$ and $d(n)-d'(n)\to0$ for four bounded sequences 
of terminating decimals implies $c(n)d(n)-c'(n)d'(n)\to0$; it follows from the rearrangement
$$
c(n)d(n)-c'(n)d'(n)=(c(n)-c'(n))d(n)+c'(n)(d(n)-d'(n)).
$$ 
The next key result gives an alternative arithmetic characterization of jumps.

\begin{pro}\label{eqcl}
Two decimals $c$ and $d$ are equivalent, $c\sim d$, if and only if $c|n-d|n\to0$.
\end{pro}
\proof
Suppose that $c\sim d$. If $c=d$ then for every $n$ even $c|n-d|n=0$. If $c<d$ or $d<c$ is a jump then 
$c|n-d|n=\pm10^{-n}$ for large $n$ and so $c|n-d|n\to0$. Suppose that $c\not\sim d$. We assume that $d<c$ is not a jump 
(the other case $c<d$ is similar). So $d<e<c$ for a decimal $e$. We have that $e|n-d|n\ge10^{-n}$ for every large $n$
and the same holds for $c|n-e|n$. Thus there is an $m\in\N$ such that $c|m-d|m\ge 2\cdot10^{-m}$. But then for every
$j>m$ we have that
$$
c|j-d|j\ge2\cdot10^{-m}-9\cdot10^{-m-1}-9\cdot10^{-m-2}-\cdots-9\cdot10^{-j}>10^{-m},
$$
which means that $(c|j-d|j)|m\ne0$. Thus $c|n-d|n\not\to 0$.
\eproof

\noindent
Decimals in a jump are distinct but infinitesimally close; this phenomenon occurs here without any 
nonstandard analysis.

We say that $c\in\D$ is a {\em hybrid limit} of a sequence $c(n)$ of terminating decimals, written $\hl_n\;c(n)=c$, if 
$$
c(n)-c|n\to 0.
$$ 
In other words, $\hl_n\;c(n)=c$ means that for every $k\in\N$ there is an $N\in\N$ such that $|c(n)-c|n|<10^{-k}$
whenever $n>N$. Hybrid limits in general are not unique: each of the three sequences $c(n)=1$, 
$c(n)=0.99\dots 9=1-10^{-n}$, and $c(1)=1$, $c(2)=0.9$, $c(3)=1$, $c(4)=0.99$, $c(5)=1$, $\dots$ 
has $\hl_n\;c(n)=1$ and $\hl_n\;c(n)=0.999\dots$; in the formal sense the first sequence converges to $1$, the second 
to $0.999\dots$ and the third is not formally convergent. To make it unique, we will regard $\hl_n\;c(n)$
as the set of all its (at most two) values. By Proposition~\ref{eqcl}, $\hl_n\;c(n)=c$ and $\hl_n\;c(n)=d$ implies
$c\sim d$, and $\hl_n\;c(n)=c$ and $c\sim d$ implies $\hl_n\;c(n)=d$. So $\hl_n\;c(n)$, if it exists, is a unique 
element of $\R$: 
$$
\hl_n\;c(n)=\emptyset\ \mbox{ or }\ \hl_n\;c(n)=\{c,d\}=[c]=[d]\in\R
$$
(here $c=d$ or $c\sim d$). Note that if $c(n),d(n)\in\T_0$ with $c(n)-d(n)\to0$ then 
$\hl_n\;c(n)=\hl_n\;d(n)$. We are ready to define addition and multiplication of real numbers.
\begin{ddef}
For $\alpha,\beta\in\R$ we set $\alpha+\beta=\hl_n\;(c|n+d|n)$ and $\alpha\beta=\hl_n\;(c|n\cdot d|n)$, 
where $c\in\alpha$ and $d\in\beta$ are arbitrary decimals.
\end{ddef}

\noindent
By the previous remarks, the results are independent of the selection of $c$ and $d$, and the resulting limits, 
if they exist, are real numbers. This completes the statement of Theorem~\ref{oR} and it remains to prove 
that the defined arithmetic operations are always defined and have required properties.

\medskip\noindent
{\bf Properties of formal and hybrid limits.} To accomplish it we first derive
a few properties of $\hl_n$. It helps to use also formal limits. We say that $c\in\D$ is a 
{\em formal limit} of a sequence of terminating decimals $c(n)$, written $\fl_n\;c(n)=c$, if for every $k\in\N$ 
there is an $N\in\N$ such that for every $n>N$ the decimals $c$ and $c(n)$ have same sign or $c$ is zero, and
$$
c(n)|k=c|k.
$$  
Formal limit, if it exists, is unique. Formal convergence is stronger than hybrid: 
$$
c=\fl_n\;c(n)\Rightarrow[c]=\hl_n\;c(n)
$$
and we presented sequence with hybrid limit but without formal limit. Further, $\fl_n\;c|n=c$ and 
$\hl_n\;c|n=[c]$ for every $c\in\D$.

\begin{pro}\label{monotonie}
If $c(n)$ and $d(n)$ are sequences of terminating decimals that have hybrid limits and $c(n)\le d(n)$ for every 
large $n$, then $$\hl_n\;c(n)\le\hl_n\;d(n).$$
\end{pro}
\proof
We denote $[c]=\hl_n\;c(n)$, $[d]=\hl_n\;d(n)$ and assume for contrary that $[c]>[d]$. Thus $c>d$ and, since 
$c\not\sim d$, there exist $k,N\in\N$ such that $c|n-d|n>10^{-k}$ for $n>N$. As $c(n)-c|n\to0$ and $d(n)-d|n\to0$, we see 
that there is an $N'$ such that $c(n)-d(n)>10^{-k-1}$ for $n>N'$. This contradicts the assumption on $c(n)$ and $d(n)$. 
\eproof

Recall that a sequence $d(n)$ of terminating decimals is {\em Cauchy} if for every $k\in\N$ there is an $N\in\N$ such that 
$|d(m)-d(n)|<10^{-k}$ whenever $m,n>N$.
\begin{pro}\label{krkonv}
Let $d(n)$ be terminating decimals.
\begin{enumerate}
\item If the sequence $d(n)$ is monotone and bounded, then it has a formal limit.
\item If the sequence $d(n)$ is Cauchy, then it has a hybrid limit.
\end{enumerate}
\end{pro}
\proof
1. We suppose that $d(n)$ is nondecreasing, the other case is similar. Using Proposition~\ref{tlup}, we set
$d=\sup_n d(n)$. We prove that $\fl_n\;d(n)=d$. If $d$ is the larger element of a jump (i.e., $d>0,d\in\T_0$ or 
$d<0,d\in\T_9$), then in fact 
$d=\max_n d(n)$ and $d(n)$ from an index on equals $d$.
Else for any given $l\in\Z$ there is a decimal $e$ such that $e<d$ and $e_i=d_i$ for $i\ge l$. Then 
$e<d(N)\le d$ for some $N$ and it follows that $d(N)_i=e_i=d_i$ for $i\ge l$. Hence, by monotonicity, the same is
true for every $d(n)$ with $n\ge N$. Thus $d(n)$ formally converge to $d$.

2. The proof is the same as in the classic case for metric limit and real numbers. Sequence $d(n)$ is Cauchy and 
therefore bounded. It contains a monotone subsequence $d(n_k)$ (every infinite sequence in a linear order has an infinite monotone subsequence). By part 1, $\fl_k\; d(n_k)=c$ for some $c\in\D$. Thus $\hl_k\; d(n_k)=[c]$. By the Cauchy property of $d(n)$ and $\hl_k\; d(n_k)=[c]$, for given $l\in\N$ there is an $N\in\N$ 
such that $m\ge n>N$ implies $|d(m)-d(n)|<10^{-l}$ and $k>N$ implies $|d(n_k)-c|k|<10^{-l}$. Thus for every $k>N$ 
we have that $|d(k)-c|k|\le|d(k)-d(n_k)|+|d(n_k)-c|k|<2\cdot10^{-k}$. Thus $d(k)-c|k\to0$ and $\hl_k\; d(k)=[c]$. 
\eproof

\begin{pro}
The sum $\alpha+\beta$ and product $\alpha\beta$ is defined for every pair of real numbers $\alpha,\beta\in\R$.
\end{pro}
\proof
If $c(n)$ and $d(n)$ are two Cauchy sequences of terminating decimals, then $c(n)+d(n)$ and $c(n)\cdot d(n)$ are 
Cauchy as well (to see the latter use the above rearrangement). For any decimals $c\in\alpha$ and $d\in\beta$ the
sequences $c|n$ and $d|n$ are Cauchy, in fact 
$$
|\,c|m-c|n\,|<10^{-n}\ \mbox{ for }\ m\ge n
$$
and similarly for $d|n$. Thus $c|n+d|n$ and $c|n\cdot d|n$ are Cauchy, have hybrid limits by 2 of 
Proposition~\ref{krkonv} and define the sum $\alpha+\beta=[c]+[d]$ and product $\alpha\cdot\beta=[c]\cdot[d]$.
\eproof

\noindent
Note that for terminating decimals $c, d$ one has $[c]+[d]=[c+d]$ and $[c][d]=[cd]$. For every $k\in\Z$ and $d\in\D$ 
one has $[10^k][d]=[10^kd]$.

We mentioned in the Introduction that it is necessary to justify exchanges of limits with sums 
and products; for hybrid limits it is easy. 
\begin{pro}\label{vymeny}
Let $c(n)$ and $d(n)$ be sequences of terminating decimals that have hybrid limits. Then so have the sequences
$c(n)+d(n)$ and $c(n)\cdot d(n)$ and
\begin{eqnarray*}
\hl_n\;(c(n)+d(n))&=&\hl_n\;c(n)+\hl_n\;d(n)\\ 
\hl_n\;(c(n)\cdot d(n))&=&\hl_n\;c(n)\cdot \hl_n\;d(n).
\end{eqnarray*}
\end{pro}
\proof
Let $\hl_n\;c(n)=[c]$, $\hl_n\;d(n)=[d]$ and $[c]+[d]=[e]$. Then $c(n)-c|n\to0$ and $d(n)-d|n\to0$, 
which sums to $c(n)+d(n)-(c|n+d|n)\to0$. By the definition of addition, $c|n+d|n-e|n\to0$. Summing again we get
$c(n)+d(n)-e|n\to0$. Thus $\hl_n\;(c(n)+d(n))=[e]$.

As for the product, let $[c]\cdot[d]=[e]$. We noted above that $c(n)-c|n\to0$ and $d(n)-d|n\to0$ implies
$c(n)\cdot d(n)-c|n\cdot d|n\to0$. By the definition of multiplication, $c|n\cdot d|n-e|n\to0$. Summing we get
$c(n)\cdot d(n)-e|n\to0$. Thus $\hl_n\;(c(n)\cdot d(n))=[e]$.
\eproof

\noindent
{\bf Conclusion of the proof of Theorem~\ref{oR}.} We need to prove that addition and multiplication on $\R$ are 
commutative and associative operations with distinct neutral elements $0$ and $1$, have inverse elements (except $0$ 
for multiplication), satisfy the distributive law, and that $(\R,<)$ is a complete linear order to which 
addition of an element and multiplication by a positive element are increasing functions.

We first prove that $(\R,+,\cdot)$ is an integral domain. It is clear that $[0]=\{0\}$ and $[1]=\{0.999\dots,1\}$ 
are neutral elements to $+$ and $\cdot$ and that the commutativity of $+$ and $\cdot$ is transferred by 
hybrid limits from $\T_0$ to $\R$. Also, it is clear that every $[c]\in\R$ has $[-c]$ as its additive inverse.
However, associativity and distributivity need Proposition~\ref{vymeny}. For example, to prove that addition is 
associative, we take $[c]$, $[d]$ and $[e]$ in $\R$ and, using Propositions~\ref{aotd} and \ref{vymeny}, get that 
indeed
\begin{eqnarray*} 
([c]+[d])+[e]&=&\hl_n\;(c|n+d|n)+\hl_n\;e|n\ \mbox{ (definition of addition)}\\
&=&\hl_n\;((c|n+d|n)+e|n)\ \mbox{ (Proposition~\ref{vymeny})}\\
&=&\hl_n\;(c|n+(d|n+e|n))\ \mbox{ (Proposition~\ref{aotd})}\\
&=&\hl_n\;c|n+\hl_n\;(d|n+e|n)\ \mbox{ (Proposition~\ref{vymeny})}\\
&=&[c]+([d]+[e])\ \mbox{ (definition of addition)}.
\end{eqnarray*}
In the same way it follows by Propositions~\ref{aotd} and \ref{vymeny} that multiplication is 
associative and that the distributive law holds. If $[c][d]=[0]$, 
then $c|n\cdot d|n\to0$, which implies (as both $c|n$ and $d|n$ is monotone) that $c|n\to0$ or $d|n\to0$, that is,
$[c]=[0]$ or $[d]=[0]$. Thus $(\R,+,\cdot)$ is an integral domain.

We consider the properties of $<$. In Proposition~\ref{tlup} we have already proven that $(\R,<)$ is a complete linear 
order. We need to show that 
$$
[c]<[d]\Rightarrow [c]+[e]<[d]+[e]\ \mbox{ and, for $[e]>0$, }\ [c][e]<[d][e].
$$ 
We have $c<d$ and therefore $c|n\le d|n$ for every $n$ 
and thus (by Proposition~\ref{aotd}) $c|n+e|n\le d|n+e|n$ for every $n$ and any real number $[e]$. Taking hybrid limits, by 
Proposition~\ref{monotonie} we get that $[c]+[e]\le[d]+[e]$. But the equality $[c]+[e]=[d]+[e]$ is impossible as 
cancelling $[e]$ (i.e., adding $[-e]$ to both sides and applying associativity of $+$) would give $[c]=[d]$. Hence 
$[c]+[e]<[d]+[e]$. For multiplication we argue similarly; the cancellation of $[e]$ in $[c][e]=[d][e]$ follows
from the fact that $(\R,+,\cdot)$ is an integral domain. Thus $(\R,+,\cdot,<)$ is an ordered ring.

It remains to show that multiplication has inverses.

\begin{pro}\label{division}
For every nonzero decimal $c$ there is a decimal $d$ such that $$[c]\cdot[d]=[1]=\{0.999\dots,1\}.$$
\end{pro}
\proof
It suffices to construct a Cauchy sequence of terminating decimals $d(n)$ with the property that 
$c|n\cdot d(n)-1\to0$. Setting $[d]=\hl_n\;d(n)$ and taking hybrid limits, we get by Proposition~\ref{vymeny}
that
$$
[0]=\hl_n\;(c|n\cdot d(n)-1)=\hl_n\;c|n\cdot\hl_n\;d(n)-\hl_n\;1=[c]\cdot[d]-[1]
$$
and (by the already established properties of $+$) $[c]\cdot[d]=[1]$.

We may assume that $c$ is positive. For given $n\in\N$, we consider decimals $c'=10^n\cdot c|n$ and 
$10^N$ in $\Z_0$, where $N\in\N$ is selected so that $10^{N-n}>c'$. Dividing in $\Z_0$, we have $10^N=c'd+e$ with 
$d,e\in\Z_0$, $d\ge0$ and 
$0\le e<c'$. Multiplying by $10^{-N}$, we get in $\T_0$ the equality
\begin{eqnarray*}
1&=&(10^{-n}c')(10^{-N+n}d)+10^{-N}e\\
&=&c|n\cdot(10^{-N+n}d)+10^{-N}e\\
&=:&c|n\cdot d(n)+e(n).
\end{eqnarray*}
As 
$$
0\le e(n)=10^{-N}e<10^{-N}c'<10^{-n}, 
$$
we have $1-c|n\cdot d(n)=e(n)\to0$. We show that $d(n)$ is Cauchy. Because
$c>0$, there exist $k\in\Z$ and $N\in\N$ such that 
$$
c|n\ge10^k\ \mbox{ for }\ n>N. 
$$
From $0<c|n\cdot d(n)=1-e(n)\le1$ we get 
$$
0<d(n)\le10^{-k}\ \mbox{ for }\ n>N. 
$$
We take $m\ge n>N$ and write $c|m=c|n+\delta$ and $d(m)=d(n)+\Delta$. Recall that
$$
0\le\delta<10^{-n}.
$$
Subtracting $1=c|n\cdot d(n)+e(n)$ from $1=(c|n+\delta)(d(n)+\Delta)+e(m)$ and rearranging, we get
\begin{eqnarray*}
10^k|\Delta|\le|\,c|m\cdot\Delta\,|&=&|e(n)-e(m)-\delta\cdot d(n)|\\
&\le&|e(n)|+|e(m)|+|\delta\cdot d(n)|\\
&<&2\cdot10^{-n}+10^{-n-k}.
\end{eqnarray*}
Thus $|\Delta|=|d(m)-d(n)|<2\cdot10^{-n-k}+10^{-n-2k}$ for $m\ge n>N$ and we see that $d(n)$ is Cauchy.
\eproof

\noindent
This concludes the proof of Theorem~\ref{oR}.

\medskip\noindent
{\bf The proof of Theorem~\ref{oQ}.} We start with an arithmetic characterization of ultimately periodic real numbers. 
For $a,b\in\N$, we let $9^{(a)}0^{(b)}$ denote the decimal $99\dots900\dots0\in\Z_0$ with $a$ $9$'s and $b$ $0$'s and 
similarly for $9^{(a)}$. 
\begin{pro}\label{archarP}
A real number $\alpha\in\R$ is in $[\P]$ if and only if $[9^{(a)}0^{(a)}]\alpha$ is in $[\Z_0]$ for some
$a\in\N$.
\end{pro}
\proof
Suppose that $\alpha=[c]$ for $c\in\P$ with period $p$. Then there is an $m\in\Z$ such that for every $k\in\N$ 
we have $[9^{(kp)}][c]=[10^{kp}][c]-[c]=[10^{kp}c]-[c]=[d]$ where $d_i=0$ for $i<m$. Taking $k$ large enough so that 
$a=kp>|m|$, we get $[10^a][9^{(a)}][c]=[9^{(a)}0^{(a)}][c]$ in $[\Z_0]$.

Now suppose that for $a\in\N$ and $c\in\D$ we have $[9^{(a)}0^{(a)}][c]\in[\Z_0]$; we may assume that $c$ is positive. 
It follows that there is a $k\in\N$ and a terminating decimal $e$ such that $[10^ac]-[c]=[e]$ and $e_i=0$ for $i<-k$.
For $l>k$ we write $(10^ac)|l=(10^ac)|k+\delta(l)$ and $c|l=c|k+\Delta(l)$. It is clear that 
$\delta(l),\Delta(l)\in\T_0$ and are 
nonnegative, smaller than $10^{-k}$, and, since $10^ac\not\in\T_0$ or $c\not\in\T_9$, there is an $m>k$ such that 
for every $l>m$ is $\delta(l)\ge10^{-m}$ or for every $l>m$ is $\Delta(l)\le10^{-k}-10^{-m}$. Thus for every $l>m$,
$$
|\,(10^ac)|l-c|l-((10^ac)|k-c|k)\,|=|\delta(l)-\Delta(l)|<10^{-k}-10^{-m}.
$$
If $(10^ac)|k-c|k\ne e$, then $|\,(10^ac)|k-c|k-e\,|\ge10^{-k}$ and, using that $e-((10^ac)|l-c|l)\to0$, by summing we get 
for large $l$ a contradiction with the displayed inequality. Hence $(10^ac)|k-c|k=e$. Since we can increase $k$, we get 
$(10^ac)|k'-c|k'=e$ for any $k'\ge k$. So $(10^ac)_i=c_i$ for every $i<-k$, which means that $c\in\P$ with period $a$.
\eproof

\begin{lem}\label{eulerfermat}
For every nonzero number $r$ in $\Z_0$ there is an $a\in\N$ such that $9^{(a)}0^{(a)}$ is divisible by $r$.
\end{lem}
\proof
We split $r=st$ where $s$ has only prime factors $2$ and $5$ and $t$ is coprime with $10=2\cdot 5$. By Euler's generalization of the little theorem of Fermat, for $p=\varphi(|t|)$ is $10^{kp}-1=9^{(kp)}$ divisible by $t$ 
for any $k\in\N$. Taking $k$ large enough so that $10^{kp}$ is divisible by $s$, we get that 
$10^{kp}(10^{kp}-1)=9^{(kp)}0^{(kp)}$ is divisible by $st=r$.
\eproof

We show that $[\P]$ is the prime field of $\R$. The prime field is formed by the elements $a/b$, where $b\ne0$ and
$a$ and $b$ are of the form  $\pm([1]+[1]+\dots+[1])$. It is clear that these are exactly the elements in $[\Z_0]$.
If $\alpha\in[\P]$, then by Proposition~\ref{archarP} $[9^{(a)}0^{(a)}]\alpha=\beta$ for some $a\in\N$ and 
$\beta\in[\Z_0]$. Thus 
$$
\alpha=\beta/[9^{(a)}0^{(a)}]
$$ 
is in the prime field. So $[\P]$ is contained in the prime field. 

Now it suffices to show that $[\P]$ is a subfield of $\R$. Let $\alpha,\beta\in[\P]$ be any 
elements. Using Proposition~\ref{archarP} and Lemma~\ref{eulerfermat}, we take numbers $a,b,c$ in $\N$ so that 
$[9^{(a)}0^{(a)}]\alpha\in[\Z_0]$, $[9^{(b)}0^{(b)}]\beta\in[\Z_0]$ and $[9^{(c)}0^{(c)}]$ is divisible by the 
product $[9^{(a)}0^{(a)}][9^{(b)}0^{(b)}]$: $[9^{(c)}0^{(c)}]=[9^{(a)}0^{(a)}][9^{(b)}0^{(b)}]r$ for 
some $r\in[\Z_0]$. Then
$$
[9^{(c)}0^{(c)}](\alpha\pm\beta)=(r[9^{(b)}0^{(b)}])([9^{(a)}0^{(a)}]\alpha)\pm(r[9^{(a)}0^{(a)}])([9^{(b)}0^{(b)}]\beta)
$$ 
is in $[\Z_0]$ and so is 
$$
[9^{(c)}0^{(c)}](\alpha\beta)=r([9^{(a)}0^{(a)}]\alpha)([9^{(b)}0^{(b)}]\beta).
$$
Thus, by Proposition~\ref{archarP}, $\alpha\pm\beta$ and $\alpha\beta$ are in $[\P]$. If $\alpha\ne0$, we take a 
$d\in\N$ so that 
$[9^{(d)}0^{(d)}]$ is in $[\Z_0]$ divisible by $[9^{(a)}0^{(a)}]\alpha$ and get that
$$
[9^{(d)}0^{(d)}](1/\alpha)=[9^{(d)}0^{(d)}][9^{(a)}0^{(a)}]/[9^{(a)}0^{(a)}]\alpha
$$
is in $[\Z_0]$, which means that $1/\alpha\in[\P]$. Thus $[\P]$ is a subfield of $\R$ and must be equal to the prime
field. This concludes the proof of Theorem~\ref{oQ}.

\medskip\noindent
{\bf Addition and multiplication of decimals.} We showed how to add and multiply real numbers $[c]$ and $[d]$ without worrying what is $c+d$ and $cd$. For completeness, we prove that formal limits define addition and multiplication 
on $\D$ correctly. We already established this with the exception of subtraction. Indeed, for any 
$c,d\in\D$ the sequence $c|n\cdot d|n$ is monotone and bounded and so $\fl_n(c|n\cdot d|n)$ exists by part 1 
of Proposition~\ref{krkonv} and defines correctly the product $cd$. This works for $c+d$ defined by $\fl_n(c|n+d|n)$
as well unless $c$ and $d$ 
have opposite signs. If this happens, $c|n+d|n$ may not be eventually monotone and we need a finer criterion of 
formal convergence.

\begin{pro}\label{finerkr}
Let $d(n)$ be a sequence of terminating decimals. Suppose that for every $k\in\N$ there exist an $N\in\N$ and 
a pair $c(k)<c(k)'$ of terminating decimals such that $c(k)_i=c(k)'_i=0$ for $i<-k$, $c(k)'-c(k)=10^{-k}$ and 
$c(k)<d(n)<c(k)'$ for every $n>N$. The sequence $d(n)$ then formally converges.
\end{pro}
\proof
Observe that 
$$
c(1)\le c(2)\le\dots\le c\le c'\le\dots\le c(2)'\le c(1)'
$$
where $c=\sup_k c(k)$ and $c'=\inf_k c(k)'$. By part 1 in Proposition~\ref{krkonv}, 
$$
\fl_k\;c(k)=c\ \mbox{ and }\ \fl_k\;c(k)'=c'. 
$$
Suppose first that $c=c'$. For given $l\in\N$ there is a $k$ such that 
$c(k)_i=c(k)'_i=c_i=c'_i$ for every $i\ge-l$. For this $k$ we have $c(k)<d(n)<c(k)'$ for every $n>N$. It follows 
that $d(n)_i=c_i=c'_i$ for every $i\ge-l$ and $n>N$. Thus $\fl_n\;d(n)=c=c'$.

The remaining case is that $c<c'$ but $c\sim c'$. We assume that $c$ and $c'$ are negative and $c\in\T_0$, $c'\in\T_9$, 
the other case is similar. It follows that $c(k)$ is eventually constant, equal to $c$. For given $l\in\N$ we fix 
a $k$ such that $c(k)=c$ and $c(k)'_i=c'_i$ for every $i\ge-l$. For this $k$ we have $c<d(n)<c(k)'$, hence 
$c'\le d(n)<c(k)'$, for every $n>N$. It follows that $d(n)_i=c'_i$ for every $i\ge-l$ and $n>N$. Thus $\fl_n\;d(n)=c'$.
\eproof

\noindent
Note that the criterion becomes invalid when the condition  $c(k)<d(n)<c(k)'$ is weakened to $c(k)\le d(n)\le c(k)'$.

\begin{pro}
Let $c,d\in\D$ be decimals. The formal limits
$$
\fl_n\;(c|n+d|n)\ \mbox{ and }\ \fl_n\;(c|n\cdot d|n)
$$
exist and define the sum and product of decimals $c+d$ and $c\cdot d$.
\end{pro}
\proof
As we mentioned, it only remains to deal with subtraction $c-d$ of positive decimals 
$c$ and $d$. We may assume that $c_i\ne d_i$ for infinitely many $i$ because else $c|n-d|n$ is eventually constant and $\fl_n\;(c|n-d|n)$ exists trivially.
Let $k\in\N$ be given. For $l>k$ we write $c|l=c|k+\delta(l)$ and $d|l=d|k+\Delta(l)$. The terminating decimals 
$\delta(l)$ and $\Delta(l)$ are nonnegative and smaller than $10^{-k}$. By the assumption on $c$ and $d$, there is an 
$l'>k$ such that $\delta(l')-\Delta(l')\ne0$. We assume that $\delta(l')-\Delta(l')<0$, the other case
is similar. It follows that $\delta(l')-\Delta(l')\le-10^{-l'}$ and $\delta(l)-\Delta(l)<0$ for every $l>l'$. 
Setting $c(k)=c|k-d|k-10^{-k}$, $c(k)'=c|k-d|k$ and $N=l'$, for $n>N$ we have 
$$
c(k)<c|n-d|n=c|k-d|k+\delta(n)-\Delta(n)<c(k)'.
$$
Thus $\fl_n\;(c|n-d|n)$ exists by Proposition~\ref{finerkr}.
\eproof

\noindent
By properties of formal and hybrid limits we have equalities 
$$
[c]+[d]=[c+d]\ \mbox{ and }\ [c][d]=[cd]
$$ 
for every pair of decimals $c,d\in\D$, not just for terminating ones. Addition and multiplication of decimals is 
trivially commutative and by Theorem~\ref{oR} also associative and distributive up to equivalence, that is to say, 
$(c+d)+e\sim c+(d+e)$, $(cd)e\sim c(de)$ and $c(d+e)\sim cd+ce$ for every $c,d,e\in\D$. In fact, for multiplication 
it holds more strongly that
$$
(cd)e=c(de)\ \mbox{ for every }c,d,e\in\D.
$$
This can be established by a version of Proposition~\ref{vymeny} for formal limits---details are left for 
the interested reader.

\medskip\noindent
{\bf References and remarks. }The oldest model of real numbers uses cuts on the
set of rational numbers and is due to Dedekind who conceived it, by his own words, on November 24, 1858 
(Dedekind \cite{dede} where he after 14 years outlines his approach). Modern version is, e.g., in Pugh 
\cite[Chapter 1.2]{pugh} or 
Rudin \cite[p. 17]{rudi}. The first rigorous theory of irrational numbers to appear in print was the construction 
based on (equivalence classes of) Cauchy sequences of rationals, due to M\'eray \cite{mera} in 1869 and later
Cantor \cite{cant} and Heine \cite{hein}. For modern treatment see, e.g., Tao \cite[Chapter 5]{tao_ma}.

Another model of $\R$ developed Weierstrass in his lectures (starting 1865/66). He published nothing on his theory 
but it was disseminated through his students, some of which (Biermann, Hettner, Hurwitz, Killing, Kossak, Pasch, 
Pincherle---see \cite{snow} and \cite[p. 46]{stol_gmei}) produced 
written accounts on his approach; thorough treatment (based on Weierstrass' courses
in 1872 and 1884) is due to Dantscher \cite{dant}. Similar theory of 
irrational numbers was developed by Stolz \cite{stol}, \cite{stol_gmei}. This approach is often dubbed as 
Weierstrass or Weierstrass--Stolz model of $\R$ based on decimal expansions (Gamelin \cite{game}, 
Kudryavtsev \cite{kudr} and elsewhere). We looked in a few primary sources
(Stolz \cite{stol} and Stolz and Gmeinder \cite{stol_gmei}, Dantscher \cite{dant} was available 
to us only in the detailed contemporary review of Miller \cite{mill}, see also Snow \cite{snow}) and, 
unsurprisingly, did not find there the modern decimal model of $\R$; decimal expansion is mostly treated not as one 
object, actually infinite string of symbols or numbers, but via partial sums as a sequence of rational 
approximations, which is not very different from the M\'eray--Cantor--Heine model. 
But in at least one case some features of decimal 
arithmetics appear---Stolz and Gmeinder \cite[pp. 48--50]{stol_gmei} give, in effect, an algorithm for
calculating the sum $c+d$ of two infinite decimals, justify subtraction and prove the associativity (up to equivalence)
of addition of infinite decimals (in the form that $-3/10^n<\{(a+b)+c\}-\{a+(b+c)\}<3/10^n$ for every $n$). However,
after that they write: ``It would come out even more tedious to explain multiplication of two real numbers. 
For building the four arithmetic operations one prefers now to switch to the Cantor's or Weierstrass' theory of 
irrational numbers.'' 

We found several outlines or sketches of the decimal model of $\R$, namely Courant and Robbins \cite{cour_robb}, 
Gamelin \cite{game}, Gowers \cite[Chapter 4]{gowe}, \cite{gowe_ro2}, and \cite{gowe_dec}, Kudryavtsev \cite{kudr}, and 
Richman \cite{rich} but no really detailed and rigorous account. This was a motivation to present
our detailed decimal model of $\R$ in Section 2. 

There are other constructions of real numbers in the literature: Eudoxus real numbers (Arthan \cite{arth}, 
A'Campo \cite{acam}, Douglas et al. \cite{doug_al} and Street \cite{stre}), the construction of de Bruijn (\cite{debr}), factorization of a set of power series (Knuth and Pratt \cite{knut_prat_loss}), an algebraic approach
of Faltin et al. (\cite{falt_al}). For more on real numbers see the books of Ebbinghaus et al. \cite{ebbi_al} 
and especially (unfortunately only for those reading Slovak) Bukovsk\'y \cite{buko}.  

\section{Irrationality of $\sqrt{2}$ in the decimal model of $\R$} 

We use the decimal models of $\R$ and $\Q$ as developed in the previous section in Theorems~\ref{oR} and \ref{oQ} 
and their proofs. In particular, recall that for every $c,d\in\T_0$, $[c]=[d]$ implies $c=d$ and that 
$$
[c]+[d]=[c+d]\ \mbox{ and }\ [c][d]=[cd].
$$
The displayed equalities hold in fact even for every $c,d\in\D$ but this generality will be needed only in the proof of 
Proposition~\ref{odm2vd} and not in Proposition~\ref{odm2vr}. Proposition~\ref{odm2vr} is on the arithmetics in 
$\R=[\D]$ based on hybrid limits and Proposition~\ref{odm2vd} is on the arithmetics in $\D$ based on formal limits. 

\begin{pro}\label{odm2vr}
In the decimal model of real numbers, the equation $x^2=[2]$ has no solution in the set of ultimately periodic 
real numbers $[\P]$.
\end{pro}
\proof
We assume for contrary that some $\alpha\in[\P]$ satisfies $\alpha^2=[2]$. By Proposition~\ref{archarP} we have 
$[9^{(a)}0^{(a)}]\alpha=[e]$ for some $a\in\N$ and $e\in\Z_0$. Thus
$$
[e]^2=([9^{(a)}0^{(a)}]\alpha)^2=[9^{(a)}0^{(a)}]^2[2]\ \mbox{ and }\ 
e^2=(9^{(a)}0^{(a)})^22.
$$
As we know from Lindstrom's proof, the last equality is impossible because the last nonzero digit on 
the left is $1,4,5,6$ or $9$ but the last nonzero digit on the right is $2$.
\eproof

\noindent
There is nothing special about $2$, the same argument works if the last nonzero digit is a quadratic non-residue modulo $10$.
We obtain the following result.

\begin{pro}
If $c$ is a positive terminating decimal whose last nonzero digit is $2, 3, 7$ or $8$ then the equation 
$x^2=[c]$ has no solution in the set of ultimately periodic real numbers $[\P]$.
\end{pro}

But we want to get the stronger result that $x^2=1.999\dots$ has no solution in $\P$. The author once devised the following ``proof''. If $d^2=1.999\dots$ for some $d\in\P$ with period $p$, then $10^pd-d=e\in\T_0$ and 
$$
e^2=(10^pd-d)^2=d^2(10^p-1)^2=1.999\dots\cdot(9^{(p)})^2.
$$
It is easy to see that any product $ab$ with $a\in\T_9$ and $b\in\T_0\backslash\{0\}$ lies in $\T_9$ and 
therefore the right side shows that $e^2$ lies in $\T_9$. This is clearly impossible because $e\in\T_0$, 
$e^2\in\T_0$ and $\T_0$ and $\T_9$ are disjoint. Hence $x^2=1.999\dots$ has in $\D$ no ultimately periodic solution.

Unfortunately, this quick argument is fallacious because of the failure of the distributive law in $\D$. 
The second equality sign in the calculation is correctly equivalence, $(10^pd-d)^2\sim d^2(10^p-1)^2$, and 
the calculation shows only that $e^2\sim1.999\dots\cdot(9^{(p)})^2$, which is 
no (immediate) contradiction. Indeed, the same argument would prove that $x^2=3.999\dots$ has no solution 
in $\P$ either, which is in conflict with the existence of the ultimately periodic solution $x=1.999\dots\;$.

\begin{pro}\label{odm2vd}
The equation $x^2=2$ has no solution in the set of decimals $\D$ and the equation $x^2=1.999\dots$ 
has no solution in the subset of ultimately periodic decimals $\P$.
\end{pro}
\proof
We assume the contrary that some decimal $d$ satisfies $d^2=2$.
Since $\fl_n\;(d|n)^2=2$, there is an $N\in\N$ such that for $n>N$ we have $((d|n)^2)_i=0$ for $i>0$ and $((d|n)^2)_0=2$. 
We fix an $m$ larger than $N$. Then $((d|m)^2)_j>0$ for some 
$j<0$ as $2$ is not a square of any terminating decimal (because of the last nonzero digits). Since 
the sequence $(d|n)^2$ is nondecreasing, for every $n\ge m$ we have $(d|n)^2\ge(d|m)^2$ and therefore 
$((d|n)^2)_i>0$ for some $i$, $j\le i<0$. But 
$\fl_n\;(d|n)^2=2$ implies that $((d|n)^2)_i=0$ for every $i$, $j\le i<0$, and large $n$, which is a contradiction.

The unsolvability of the equation $x^2=1.999\dots$ in $\P$ follows immediately from Proposition~\ref{odm2vr}: if some 
$d\in\P$ satisfies $d^2=1.999\dots$, then $[d]^2=[d^2]=[1.999\dots]=[2]$ in $\R$, which was shown to be impossible.
\eproof

\noindent
The completeness of $(\D,<)$ provides the aperiodic solution(s) of $x^2=1.999\dots,$ 
$$
x=\pm1.414213562373095048801688724209698\dots
$$ 
For algorithms calculating this decimal see Sebah and Gourdon \cite{seba_gour}.

\end{document}